\documentclass[11pt]{article} 
\usepackage{graphicx}
\usepackage{latexsym}
\usepackage{amsmath}
\usepackage{amssymb}
\usepackage{amsthm}
\usepackage{verbatim}
\usepackage{tikz}

\title{Shrinking Circular Nim}
\author{Hiromi Oginuma\ \thanks{Nara Women's University.
{\tt xah\_oginuma@cc.nara-wu.ac.jp}} \and
Masato Shinoda \thanks{Nara Women's University. {\tt shinoda@cc.nara-wu.ac.jp}}}

\date{\today}

\begin{document}

\theoremstyle{definition} 
\newtheorem{theorem}{Theorem}[section]
\newtheorem{definition}[theorem]{Definition}
\newtheorem{lemma}[theorem]{Lemma}
\newtheorem{proposition}[theorem]{Proposition}
\newtheorem{corollary}[theorem]{Corollary}
\newtheorem{example}[theorem]{Example}
\newtheorem{remark}[theorem]{Remark}
\newtheorem{conjecture}[theorem]{Conjecture}

\maketitle              
\begin{abstract}
The game of Nim, which has been well known for many years, has numerous variations. One such variation is Circular Nim, where piles of stones are arranged on a circumference, and players take stones from consecutive adjacent piles in one move. In this paper, we propose a new variant called Shrinking Circular Nim, in which the size of the circle decreases as the game progresses. We also examine the winning conditions in this variant. This paper discusses that Shrinking Circular Nim is solved for initial configurations with five or fewer piles. In addition, we derive the winning condition for the specific case with eight piles.

\end{abstract}
\section{Introduction}

Nim is a game that has been well known for many years. In the most basic version of Nim, there are three piles of stones, and two players take turns. On their turn, a player selects one pile and removes any positive number of stones from that pile. The configuration of the game during play is represented as $(x, y, z)$, listing the number of stones in the three piles (where $x$, $y$, and $z$ are all non-negative integers). 
For example, from $(3, 5, 7)$, the player whose turn it is can move to $(3, 2, 7)$, $(0, 5, 7)$, and so on. The players take turns making these moves, and the player who takes the last stone wins the game. For instance, from the configuration $(0, 0, 3)$, the player whose turn it is can move to $(0, 0, 0)$, the terminal configuration, and win the game. This game can also be played with more than three piles using the same rules.

Nim is a two-player, zero-sum, finite, deterministic, and perfect information game with no possibility of a draw. Therefore, every configuration can be classified as either an N-position, where the player whose turn it is can force a win, or a P-position, where the player whose turn it is will inevitably lose if the opponent plays optimally. The classification of N-positions and P-positions in the aforementioned version of Nim was demonstrated by Bouton
\cite{bou02} 
in 1902. It became known that the winning conditions in this game include mathematically intriguing structures.

In this study, we propose Shrinking Circular Nim, a modification of the rules of Circular Nim introduced by Dufour-Heubach \cite{duf13}, where piles of stones are arranged on a circumference, and stones are taken from multiple adjacent piles simultaneously in one move. We examine the winning conditions in Shrinking Circular Nim.

The organization of this paper is as follows. In Chapter 2, we will first explain the rules of Shrinking Circular Nim while comparing them to those of Circular Nim. Then, in Chapter 3, we will present the results we obtained regarding the winning conditions for Shrinking Circular Nim, and in Chapter 4, we will provide the proof of the main theorems.

\section{Circular Nim and Shirinking Circular Nim}
\subsection{Circular Nim}

First, we will explain the rules of Circular Nim and describe the current state of research. In all the games discussed in this paper, there are two players who take turns removing stones, and the player who removes the last stone wins (if a player has no stones to remove on their turn, they lose). 

\begin{definition}[Circular Nim]
Let $n$ and $k$ be positive integers such that $k \leq n$. In this game, $n$ piles of stones are arranged on a circumference. Each player, on their turn, selects $k$ consecutive adjacent piles and removes any number of stones from each pile. The total number of stones removed must be at least one, but it is permissible for some of the $k$ chosen piles to have no stones removed. This game is denoted as ${\rm CN}(n, k)$.
\end{definition}

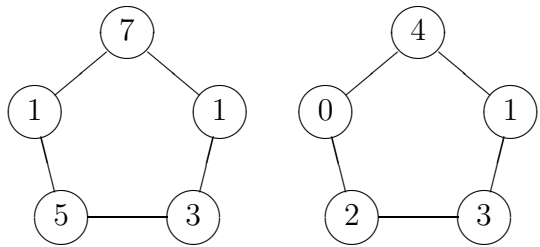
\begin{figure}[ht]
\unitlength.1pt
\begin{picture}(1100,900)(-600,-100)
\put(400,0){\circle{200}}
\put(900,0){\circle{200}}
\put(300,400){\circle{200}}
\put(1000,400){\circle{200}}
\put(650,700){\circle{200}}
\put(500,0){\line(1,0){300}}
\put(376,97){\line(-1,4){52}}
\put(924,97){\line(1,4){52}}
\put(378,465){\line(6,5){194}}
\put(922,465){\line(-6,5){194}}
\put(370,-30){\large{$5$}}
\put(870,-30){\large{$3$}}
\put(270,370){\large{$1$}}
\put(970,370){\large{$1$}}
\put(620,670){\large{$7$}}
\put(1500,0){\circle{200}}
\put(2000,0){\circle{200}}
\put(1400,400){\circle{200}}
\put(2100,400){\circle{200}}
\put(1750,700){\circle{200}}
\put(1600,0){\line(1,0){300}}
\put(1476,97){\line(-1,4){52}}
\put(2024,97){\line(1,4){52}}
\put(1478,465){\line(6,5){194}}
\put(2022,465){\line(-6,5){194}}
\put(1470,-30){\large{$2$}}
\put(1970,-30){\large{$3$}}
\put(1370,370){\large{$0$}}
\put(2070,370){\large{$1$}}
\put(1720,670){\large{$4$}}
\end{picture}
\caption{An example of game progression of ${\rm CN}(5,3)$}
\end{figure}

As a concrete example of the game, we will explain a configuration in ${\rm CN}(5,3)$. In the left diagram in Figure 1, piles of stones are arranged on a circumference (the numbers inside the small circles represent the number of stones in each pile). From this configuration, by taking stones from each of three consecutive adjacent piles, it is possible to move to the configuration shown in the right diagram above. Keep in mind that even if some piles are depleted of stones, the circle is not compressed (the original adjacency relationship among the five piles is maintained). Therefore, in the next turn for the configuration shown in the right diagram above, it is not possible to take stones simultaneously from the piles with 4, 2, and 3 stones.

Depending on the values of $n$ and $k$, Circular Nim can be a trivial game or equivalent to the classic Nim or its variations, with known conditions for determining the winner. 
In particular, ${\rm CN}(n,n-1)$ is a special case of Moore's Nim (Moore\cite{moo10}), and the P-positions take the the configuration of the form $(a,a,a,\ldots,a)$. 
The new research focuses on cases where $2\leq k \leq n-2$. 
So far, the cases of $(n, k) = (4, 2), (5, 2), (5, 3), (6, 3), (6, 4), (7, 4), (8, 6)$ have been solved
(\cite{duf13}, 
\cite{duf22}).

\subsection{Shrinking Circular Nim}

In this study, we define Shrinking Circular Nim, a game in which some of the rules of Circular Nim are modified.

\begin{definition}[Shrinking Circular Nim] 
Let $n$ and $k$ be positive integers such that $k\leq n$. In the initial configuration, $n$ piles of stones are arranged on a circumference. Each player, on their turn, selects $k$ consecutive adjacent piles and removes any number of stones from each pile (with a total of at least one stone removed). However, if any pile is reduced to zero stones during the game, that pile is considered to be vanished. The players then take stones from $k$ consecutive adjacent piles from the remaining piles of stones. This game is denoted as ${\rm SCN}(n, k)$.
\end{definition}

The main difference between Circular Nim and Shrinking Circular Nim is as follows: In Circular Nim, when a pile of stones is depleted, its position is not closed up. In contrast, in Shrinking Circular Nim, the position of a depleted pile is closed up, resulting in a smaller circle. For example, in the left diagram in Figure 2, if one stone is removed from the pile with three stones and one stone is removed from the pile with one stone, the number of piles is reduced by one, as shown in the right diagram.
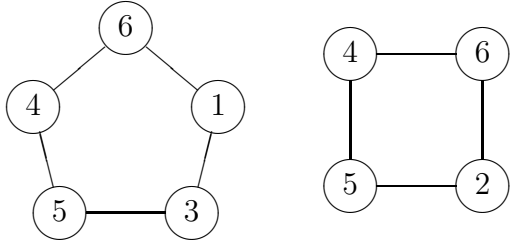
\begin{figure}[ht]
\unitlength.1pt
\begin{picture}(2100,950)(-700,-100)
\put(400,0){\circle{200}}
\put(900,0){\circle{200}}
\put(300,400){\circle{200}}
\put(1000,400){\circle{200}}
\put(650,700){\circle{200}}
\put(500,0){\line(1,0){300}}
\put(376,97){\line(-1,4){52}}
\put(924,97){\line(1,4){52}}
\put(378,465){\line(6,5){194}}
\put(922,465){\line(-6,5){194}}
\put(370,-30){\large{$5$}}
\put(870,-30){\large{$3$}}
\put(270,370){\large{$4$}}
\put(970,370){\large{$1$}}
\put(620,670){\large{$6$}}
\unitlength.1pt
\put(1500,100){\circle{200}}
\put(2000,100){\circle{200}}
\put(1500,600){\circle{200}}
\put(2000,600){\circle{200}}
\put(1600,100){\line(1,0){300}}
\put(1600,600){\line(1,0){300}}
\put(1500,200){\line(0,1){300}}
\put(2000,200){\line(0,1){300}}
\put(1470,70){\large{$5$}}
\put(1970,70){\large{$2$}}
\put(1470,570){\large{$4$}}
\put(1970,570){\large{$6$}}
\end{picture}
\caption{An example of game progression of ${\rm SCN}(5,3)$}
\end{figure}

This new rule can be considered a natural setting when actually playing this game. This is because, when playing this stone-taking game, unless the positions of the piles of stones are predetermined and fixed, a cluster of gathered stones would be recognized as a single pile. Therefore, a position where the stones have been depleted is regarded as having no pile. In this paper, we assume that the number of stones in each pile in ${\rm CN}(n,k)$ is a non-negative integer, and the number of stones in each pile in ${\rm SCN}(n,k)$ is a positive integer. Configurations that coincide through rotation or reflection (for example, 
(1,2,3,4), (2,3,4,1) and (4,3,2,1)) are considered identical.

The simplest case where the set of P-positions differs between Circular Nim and Shrinking Circular Nim is described in the following theorem.

\begin{theorem}
The P-positions in ${\rm SCN}(4,2)$ are the configurations of the forms $(a, b, a, b)$ (where $a \neq b$) and $(a, a, a)$, as well as the terminal configuration $()$. 

\begin{figure}[ht]
\unitlength.1pt
\begin{picture}(2100,900)(-650,-50)
\put(1400,200){\circle{200}}
\put(1900,200){\circle{200}}
\put(1650,700){\circle{200}}
\put(1500,200){\line(1,0){300}}
\put(1445,290){\line(1,2){160}}
\put(1855,290){\line(-1,2){160}}
\put(1370,170){\large{$a$}}
\put(1870,170){\large{$a$}}
\put(1620,670){\large{$a$}}
\put(400,200){\circle{200}}
\put(900,200){\circle{200}}
\put(400,700){\circle{200}}
\put(900,700){\circle{200}}
\put(500,200){\line(1,0){300}}
\put(500,700){\line(1,0){300}}
\put(400,300){\line(0,1){300}}
\put(900,300){\line(0,1){300}}
\put(370,170){\large{$b$}}
\put(870,170){\large{$a$}}
\put(370,670){\large{$a$}}
\put(870,670){\large{$b$}}
\put(550,0){$a\neq b$}
\end{picture}
\caption{P-positions of  ${\rm SCN}(4,2)$}
\end{figure}
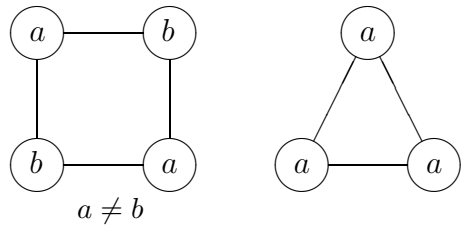
\end{theorem}

Note that the set of P-positions for ${\rm SCN}(4,2)$ also includes P-positions where the number of piles is reduced during the game (which are also P-positions for ${\rm SCN}(3,2)$). The non-terminal P-positions are shown in Figure 3.

\section{Main Theorems}
\subsection{Shrinking Circular Nim with Five Piles}

The set of P-positions when the number of piles is five can be determined as follows. Note that the cases other than $(n, k) = (5, 2)$ and $(5, 3)$ are trivial and thus omitted.

\subsubsection{${\rm SCN}(5,2)$}
In this subsection, we will discuss the case of Shrinking Circular Nim with $(n, k) = (5, 2)$. First, for reference, we will describe the set of P-positions for ${\rm CN}(5, 2)$.

\begin{theorem}[Dufour\cite{duf96},Ehrenborg-Steingr\'{i}msson\cite{ehr96}]
The P-positions in ${\rm CN}(5, 2)$ are the configurations shown in Figure 4, that is, configurations of the form $(M, m, a, b, m)$, where $M=\max\{M,a,b,m\}$ and $M+m=a+b$ are satisfied.
\end{theorem}
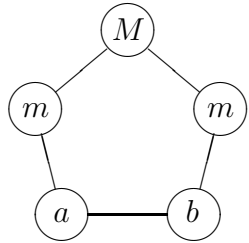
\begin{figure}[ht]
\unitlength.1pt
\begin{picture}(1100,1050)(-1150,-250)
\put(400,0){\circle{200}}
\put(900,0){\circle{200}}
\put(300,400){\circle{200}}
\put(1000,400){\circle{200}}
\put(650,700){\circle{200}}
\put(500,0){\line(1,0){300}}
\put(376,97){\line(-1,4){52}}
\put(924,97){\line(1,4){52}}
\put(378,465){\line(6,5){194}}
\put(922,465){\line(-6,5){194}}

\put(370,-30){\large{$a$}}
\put(870,-30){\large{$b$}}
\put(250,370){\large{$m$}}
\put(950,370){\large{$m$}}
\put(590,660){\large{$M$}}

\put(-250,-220){$m\leq a\leq M,m\leq b\leq M,m+M=a+b$}
\end{picture}
\caption{P-positions of  ${\rm CN}(5,2)$}
\end{figure}

In the configuration shown in Figure 4, the largest number of stones in the five piles is $M$, and the smallest number is $m$. Thus all P-positions in ${\rm CN}(5,2)$ are given by a simple formula.  In contrast, in Shrinking Circular Nim, it is found that the set of P-positions includes exceptional configurations as described below.

\begin{theorem}
The P-positions in ${\rm SCN}(5,2)$ are the two types of configurations shown in Figure 5 when  $m$ is even, and the four types of configurations shown in Figure 6 when  $m$ is odd, as well as the P-positions from ${\rm SCN}(4,2)$: the configurations of the forms $(a, b, a, b)$ (where $a \neq b$),  $(a, a, a)$ and the terminal configuration $()$. 

\begin{figure}[ht]
\unitlength.1pt
\begin{picture}(2200,1200)(-600,-350)
\put(400,0){\circle{200}}
\put(900,0){\circle{200}}
\put(300,400){\circle{200}}
\put(1000,400){\circle{200}}
\put(650,700){\circle{200}}
\put(500,0){\line(1,0){300}}
\put(376,97){\line(-1,4){52}}
\put(924,97){\line(1,4){52}}
\put(378,465){\line(6,5){194}}
\put(922,465){\line(-6,5){194}}

\put(370,-30){\large{$a$}}
\put(870,-30){\large{$b$}}
\put(250,370){\large{$m$}}
\put(950,370){\large{$m$}}
\put(590,670){\large{$M$}}
\put(50,-200){$m<a<M,m<b<M$,}
\put(250,-320){$m+M=a+b$}
\put(1540,0){\circle{200}}
\put(2000,0){\circle{200}}
\put(1400,400){\circle{200}}
\put(2100,400){\circle{200}}
\put(1750,700){\circle{200}}
\put(1630,0){\line(1,0){280}}
\put(1476,97){\line(-1,4){52}}
\put(2024,97){\line(1,4){52}}
\put(1478,465){\line(6,5){194}}
\put(2022,465){\line(-6,5){194}}
\put(1500,-30){\large{$m$}}
\put(1950,-30){\large{$m$}}
\put(1340,370){\large{$M$}}
\put(2040,370){\large{$M$}}
\put(1690,670){\large{$m+1$}}
\put(1500,-200){$m+2\leq M$}
\end{picture}
\caption{P-positions of  ${\rm SCN}(5,2)$ when $m$ is even}
\end{figure}
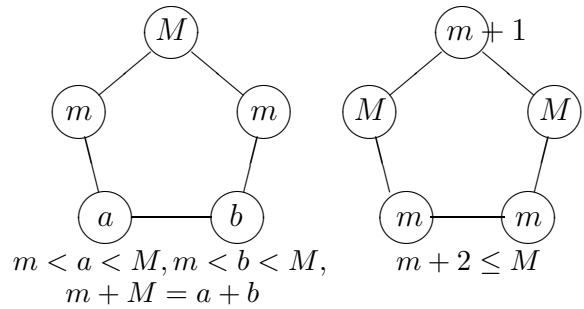

\begin{figure}[h]
\unitlength.1pt
\begin{picture}(2300,2200)(-600,-200)
\put(400,1200){\circle{200}}
\put(900,1200){\circle{200}}
\put(300,1600){\circle{200}}
\put(1000,1600){\circle{200}}
\put(650,1900){\circle{200}}
\put(500,1200){\line(1,0){300}}
\put(376,1297){\line(-1,4){52}}
\put(924,1297){\line(1,4){52}}
\put(378,1665){\line(6,5){194}}
\put(922,1665){\line(-6,5){194}}

\put(370,1170){\large{$a$}}
\put(870,1170){\large{$b$}}
\put(250,1570){\large{$m$}}
\put(950,1570){\large{$m$}}
\put(590,1870){\large{$M$}}
\put(50,1000){$m<a<M,m<b<M$,}
\put(250,880){$m+M=a+b$}
\put(1500,1200){\circle{200}}
\put(2000,1200){\circle{200}}
\put(1400,1600){\circle{200}}
\put(2100,1600){\circle{200}}
\put(1750,1900){\circle{200}}
\put(1600,1200){\line(1,0){300}}
\put(1476,1297){\line(-1,4){52}}
\put(2024,1297){\line(1,4){52}}
\put(1478,1665){\line(6,5){194}}
\put(2022,1665){\line(-6,5){194}}
\put(1450,1170){\large{$m$}}
\put(1950,1170){\large{$m$}}
\put(1340,1570){\large{$M$}}
\put(2040,1570){\large{$M$}}
\put(1690,1870){\large{$m+1$}}
\put(1500,1000){$m+3\leq M$}

\put(400,0){\circle{200}}
\put(900,0){\circle{200}}
\put(300,400){\circle{200}}
\put(1000,400){\circle{200}}
\put(650,700){\circle{200}}
\put(500,0){\line(1,0){300}}
\put(376,97){\line(-1,4){52}}
\put(924,97){\line(1,4){52}}
\put(378,465){\line(6,5){194}}
\put(922,465){\line(-6,5){194}}
\put(350,-30){\large{$m$}}
\put(850,-30){\large{$m$}}
\put(240,370){\large{$m+1$}}
\put(940,370){\large{$m+1$}}
\put(590,670){\large{$m+2$}}
\put(1500,0){\circle{200}}
\put(2000,0){\circle{200}}
\put(1400,400){\circle{200}}
\put(2100,400){\circle{200}}
\put(1750,700){\circle{200}}
\put(1600,0){\line(1,0){300}}
\put(1476,97){\line(-1,4){52}}
\put(2024,97){\line(1,4){52}}
\put(1478,465){\line(6,5){194}}
\put(2022,465){\line(-6,5){194}}
\put(1450,-30){\large{$m$}}
\put(1950,-30){\large{$m$}}
\put(1340,370){\large{$m+1$}}
\put(2040,370){\large{$m+2$}}
\put(1690,670){\large{$m+1$}}
\end{picture}
\caption{P-positions of  ${\rm SCN}(5,2)$ when $m$ is odd}
\end{figure}
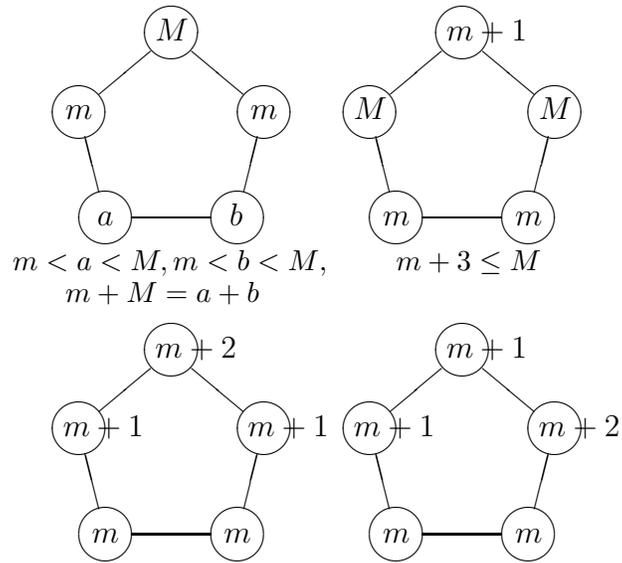

\end{theorem}

We again  remark that in each configuration in  Figure 5 and  Figure 6 
the maximum number of stones among the five piles is denoted $M$, and the minimum as $m$.

\subsubsection{${\rm SCN}(5,3)$}
In this subsection, we will discuss the case of Shrinking Circular Nim with $(n, k) = (5, 3)$. Here again,  for reference, we will describe the set of P-positions for ${\rm CN}(5, 3)$.

\begin{theorem}[\cite{duf96},\cite{ehr96}]
The P-positions in ${\rm CN}(5, 3)$ are those shown in Figure 7, that is, configurations of the form $(0,M, a, b, M)$, where $M=\max\{M,a,b\}$ and $M=a+b$ are satisfied.
\begin{figure}[ht]
\unitlength.1pt
\begin{picture}(1100,1000)(-1100,-250)
\put(400,0){\circle{200}}
\put(900,0){\circle{200}}
\put(300,400){\circle{200}}
\put(1000,400){\circle{200}}
\put(650,700){\circle{200}}
\put(500,0){\line(1,0){300}}
\put(376,97){\line(-1,4){52}}
\put(924,97){\line(1,4){52}}
\put(378,465){\line(6,5){194}}
\put(922,465){\line(-6,5){194}}

\put(370,-30){\large{$a$}}
\put(870,-30){\large{$b$}}
\put(240,370){\large{$M$}}
\put(940,370){\large{$M$}}
\put(610,660){\large{$0$}}
\put(380,-220){$M=a+b$}
\end{picture}
\caption{P-positions of  ${\rm CN}(5,3)$}
\end{figure}
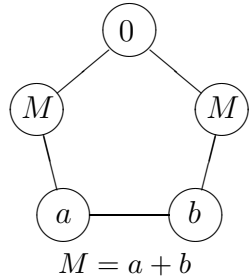

\end{theorem}


 In contrast, in Shrinking Circular Nim, it is found that the set of P-positions includes exceptional  configurations as described below.

\begin{theorem}
The P-positions in ${\rm SCN}(5,3)$ are the  configurations shown in Figure 8, as well as the P-positions in ${\rm SCN}(4,3)$: configuration of the form $(a, a, a, a)$ and the terminal configuration $()$.
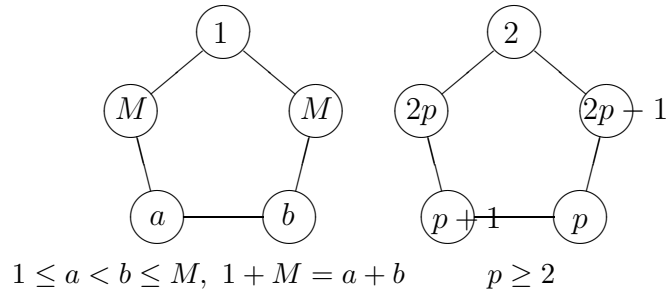
\begin{figure}[ht]
\unitlength.1pt
\begin{picture}(2200,1050)(-600,-280)
\put(400,0){\circle{200}}
\put(900,0){\circle{200}}
\put(300,400){\circle{200}}
\put(1000,400){\circle{200}}
\put(650,700){\circle{200}}
\put(500,0){\line(1,0){300}}
\put(376,97){\line(-1,4){52}}
\put(924,97){\line(1,4){52}}
\put(378,465){\line(6,5){194}}
\put(922,465){\line(-6,5){194}}
\put(370,-30){\large{$a$}}
\put(870,-30){\large{$b$}}
\put(240,370){\large{$M$}}
\put(940,370){\large{$M$}}
\put(610,660){\large{$1$}}
\put(-150,-250){$1\leq a<b\leq M, \ 1+M=a+b$}
\put(1500,0){\circle{200}}
\put(2000,0){\circle{200}}
\put(1400,400){\circle{200}}
\put(2100,400){\circle{200}}
\put(1750,700){\circle{200}}
\put(1600,0){\line(1,0){300}}
\put(1476,97){\line(-1,4){52}}
\put(2024,97){\line(1,4){52}}
\put(1478,465){\line(6,5){194}}
\put(2022,465){\line(-6,5){194}}
\put(1440,-40){\large{$p+1$}}
\put(1970,-40){\large{$p$}}
\put(1340,370){\large{$2p$}}
\put(2010,370){\large{$2p-1$}}
\put(1710,660){\large{$2$}}
\put(1650,-250){$p\geq 2$}
\end{picture}
\caption{P-positions of  ${\rm SCN}(5,3)$}
\end{figure}
\end{theorem}

\begin{remark}
If we subtract one stone from each pile in a typical P-position $(1,M,a,b,M)$ in ${\rm SCN}(5,3)$, it becomes 
$(0,M-1,a-1,b-1,M-1)$, which satisfies the P-position condition for ${\rm CN}(5,3)$, namely $(M-1)=(a-1)+(b-1)$.
This is significant because in ${\rm SCN}(5,3)$, when it is not possible to eliminate a pile to move to a 
P-position in ${\rm SCN}(4,3)$, we only consider moves that leave at least one stone in each pile. Therefore, the relationship in this modified configuration (after subtracting one from each pile) becomes important. This correspondence between Circular Nim and Shrinking Circular Nim is also used in proving the winning conditions for ${\rm SCN}(8,6)$.
\end{remark}

\subsection{Shrinking Circular Nim with Six Piles}
In this chapter we examine Shrinking Circular Nim when the number of piles is six. 
\subsubsection{${\rm SCN}(6,2)$}
First, in the case of $(n, k) = (6, 2)$, the winning conditions for Circular Nim remain unsolved. As we saw in the previous chapter, the winning conditions for Shrinking Circular Nim are considered to be more complex than those of Circular Nim, and thus the conditions for ${\rm SCN}(6,2)$ also remain unsolved.
\subsubsection{${\rm SCN}(6,3)$}
For Shrinking Circular Nim with $(n,k)=(6,3)$, it is known that the set of 
P-positions for the corresponding Circular Nim has a notable pattern: the differences in the number of stones between opposite piles on the circumference are equal, that is, configurations of the form  $(a,b+q,c,a+q,b,c+q)$ using integers $a,b,c$ and $q\geq 0$ (See \cite{duf13}). However, we have not yet determined the set of P-positions for ${\rm SCN}(6,3)$.  
It is clear that it includes all P-positions in ${\rm SCN}(5,3)$. Through computer analysis of specific configurations with all six piles remain, we observe that many configurations in the P-position set, the differences in the number of stones between each pair of opposite piles are equal. 
However, even if these differences are all equal, certain configurations are not P-positions 
(for example, $(a,b,a+b-1,a,b,a+b-1)$ is not a P-position since it can be moved to $(1,a+b-1,a,b,a+b-1)$ when $a\neq b$). 
Conversely, there are also many P-positions that do not satisfy this condition (for example, $(1,6,2,3,3,6)$ is a P-position). 

\subsubsection{${\rm SCN}(6,4)$}
In this subsection we discuss Shrinking Circular Nim in the case of $(n,k)=(6,4)$. The set of P-positions for the corresponding Circular Nim includes the condition from ${\rm CN}(6,3)$, that the differences in the number of stones between each pair of opposite piles are equal, along with an additional condition regarding the Nim-sum that must also be satisfied.

\begin{theorem}[\cite{duf13}]
The P-positions in ${\rm CN}(6, 4)$ are the  configurations shown in Figure 9, that is,  configurations of the form  $(a,b+q,c,a+q,b,c+q)$ where  $a\oplus b\oplus c=0$ are satisfied.

\begin{figure}[h]
\unitlength.1pt
\begin{picture}(1100,1100)(-1100,-250)
\put(700,0){\circle{200}}
\put(400,180){\circle{200}}
\put(1000,180){\circle{200}}
\put(400,540){\circle{200}}
\put(1000,540){\circle{200}}
\put(700,720){\circle{200}}
\put(400,280){\line(0,1){160}}
\put(1000,280){\line(0,1){160}}
\put(486,129){\line(5,-3){130}}
\put(914,129){\line(-5,-3){130}}
\put(486,591){\line(5,3){130}}
\put(914,591){\line(-5,3){130}}
\put(630,-30){\large{$a+q$}}
\put(370,150){\large{$c$}}
\put(970,150){\large{$b$}}
\put(330,510){\large{$b+q$}}
\put(930,510){\large{$c+q$}}
\put(670,690){\large{$a$}}
\put(400,-220){$a\oplus b \oplus c=0$}
\end{picture}
\caption{P-positions of  ${\rm CN}(6,4)$}
\end{figure}
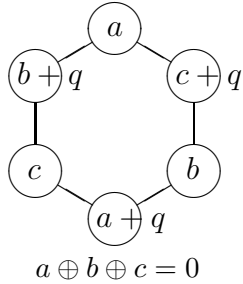
\end{theorem}

The set of P-positions for ${\rm SCN}(6,4)$ has also not been determined in this paper. Through computer analysis of specific configurations with all six piles remain, we propose the following conjecture.

\begin{conjecture}  
The P-positions in ${\rm SCN}(6,4)$ are classified into one of the following categories: (i), (ii), or (iii): 
\begin{enumerate}
\renewcommand{\theenumi}{(\roman{enumi})}
\item P-positions in ${\rm SCN}(5,4)$, that is,  configuration of the form $(a,a,a,a,a)$ the terminal configuration $()$.
\item Configurations that can be expressed as $(a,b+q,c,a+q,b,c+q)$ for some $a,b$ and $q\geq 0$, with at most one corresponding value of $c$.
\item Specific configurations $(5,9,10,7,8,12), (5,10,11,7,9,13)$ and $(5,11,11,$\\$8,9,14)$.
\end{enumerate}
\end{conjecture}

\subsection{Shrinking Circular Nim with Seven or More Piles}
In this chapter, we discuss Shrinking Circular Nim with seven or more piles, focusing specially on cases where the winning conditions for ${\rm CN}(n,k)$ have been determined, namely $(n,k)=(7,4)$ and $(8,6)$.

\subsubsection{${\rm SCN}(7,4)$}
To solve Shrinking Circular Nim for $(n,k)=(7,4)$, knowledge of the P-positions in ${\rm SCN}(6,4)$ is essential for situations where the number of piles decreases during gameplay. However, as noted in the previous chapter, the winning conditions for ${\rm SCN}(6,4)$ remain unsolved, making it challenging to determine the winning conditions for 
${\rm SCN}(7,4)$ at this stage.

\subsubsection{${\rm SCN}(8,6)$}
For Shrinking Circular Nim with $(n,k)=(8,6)$, we have been able to fully determine the conditions for identifying P-positions. To explain this, we first reference the set of P-positions for ${\rm CN}(8,6)$.

\begin{theorem}[\cite{duf13}]
The P-positions in ${\rm CN}(8,6)$ are the configurations shown in Figure 10, that is, configurations of the form $(0,M,a,M-a,\alpha,M-b,b,M)$, where 
$\alpha=\min\{M,a+b\}$ and $M$ is the maximum number of stones in any of the eight piles.
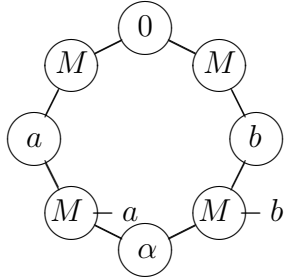
\begin{figure}[ht]
\unitlength.1pt
\begin{picture}(1100,1200)(-1200,-200)
\put(500,80){\circle{200}}
\put(220,220){\circle{200}}
\put(220,780){\circle{200}}
\put(80,500){\circle{200}}
\put(500,920){\circle{200}}
\put(780,220){\circle{200}}
\put(780,780){\circle{200}}
\put(920,500){\circle{200}}
\put(309,175){\line(2,-1){102}}
\put(691,175){\line(-2,-1){102}}
\put(309,825){\line(2,1){102}}
\put(691,825){\line(-2,1){102}}
\put(175,309){\line(-1,2){54}}
\put(825,309){\line(1,2){54}}
\put(175,691){\line(-1,-2){54}}
\put(825,691){\line(1,-2){54}}
\put(470,890){\large{$0$}}
\put(470,50){\large{$\alpha$}}
\put(140,190){\large{$M-a$}}
\put(700,190){\large{$M-b$}}
\put(50,470){\large{$a$}}
\put(890,470){\large{$b$}}
\put(160,750){\large{$M$}}
\put(720,750){\large{$M$}}
\put(-200,-160){$\alpha=\min\{M,a+b\}$, $a,b\leq M$}
\end{picture}
\caption{P-positions of  ${\rm CN}(8,6)$}
\end{figure}
\end{theorem}

In contrast, the set of P-positions for ${\rm SCN}(8,6)$ is as follows:

\begin{theorem}
The P-positions in ${\rm SCN}(8,6)$ are the  configurations shown in the left diagram of Figure 11; specifically, these are  configurations of the form $(1,M,a,M-a+1,\alpha,M-b+1,b,M)$ where $a,b\leq M$ and $\alpha=\min\{M,a+b-1\}$, as well as the P-positions in ${\rm SCN}(7,6)$: configuration of the form $(a,a,a,a,a,a,a)$ and the terminal configuration $()$. However, configurations of the form $(1,2p-1,p,p,2p-1,p,p,2p-1)$ (where $p\geq 1$, as shown in the right diagram of  Figure 11) are excluded.
\end{theorem}
\begin{figure}[ht]
\unitlength.1pt
\begin{picture}(2200,1200)(-800,-200)
\put(500,80){\circle{200}}
\put(220,220){\circle{200}}
\put(220,780){\circle{200}}
\put(80,500){\circle{200}}
\put(500,920){\circle{200}}
\put(780,220){\circle{200}}
\put(780,780){\circle{200}}
\put(920,500){\circle{200}}
\put(309,175){\line(2,-1){102}}
\put(691,175){\line(-2,-1){102}}
\put(309,825){\line(2,1){102}}
\put(691,825){\line(-2,1){102}}
\put(175,309){\line(-1,2){54}}
\put(825,309){\line(1,2){54}}
\put(175,691){\line(-1,-2){54}}
\put(825,691){\line(1,-2){54}}
\put(470,890){\large{$1$}}
\put(470,50){\large{$\alpha$}}
\put(-30,190){\large{$M-a+1$}}
\put(530,190){\large{$M-b+1$}}
\put(50,470){\large{$a$}}
\put(890,470){\large{$b$}}
\put(160,750){\large{$M$}}
\put(720,750){\large{$M$}}
\put(-250,-140){$a,b\leq M$, $\alpha=\min\{M,a+b-1\}$}
\put(1700,80){\circle{200}}
\put(1420,220){\circle{200}}
\put(1420,780){\circle{200}}
\put(1280,500){\circle{200}}
\put(1700,920){\circle{200}}
\put(1980,220){\circle{200}}
\put(1980,780){\circle{200}}
\put(2120,500){\circle{200}}
\put(1509,175){\line(2,-1){102}}
\put(1891,175){\line(-2,-1){102}}
\put(1509,825){\line(2,1){102}}
\put(1891,825){\line(-2,1){102}}
\put(1375,309){\line(-1,2){54}}
\put(2025,309){\line(1,2){54}}
\put(1375,691){\line(-1,-2){54}}
\put(2025,691){\line(1,-2){54}}
\put(1670,890){\large{$1$}}
\put(1650,50){\large{$2p-1$}}
\put(1390,190){\large{$p$}}
\put(1950,190){\large{$p$}}
\put(1250,470){\large{$p$}}
\put(2090,470){\large{$p$}}
\put(1370,750){\large{$2p-1$}}
\put(1920,750){\large{$2p-1$}}
\put(1600,-140){$p\geq 1$}
\end{picture}
\caption{P-positions of  ${\rm SCN}(8,6)$ and exception}
\end{figure}
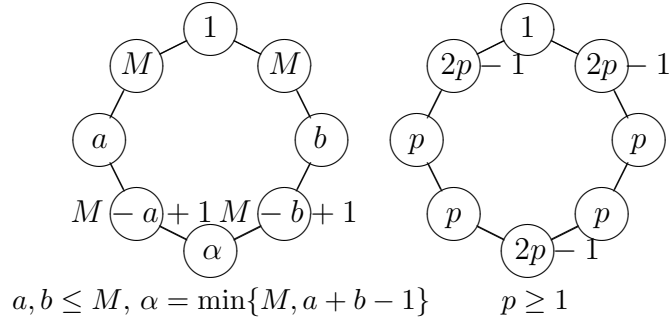

\begin{remark}
The configuration in the right diagram of Figure 11 represents a special case of the configuration in the left diagram, where $M=2p-1$ and $a=b=p.$
\end{remark}

\section{Proofs of the Main Theorems}
In this chapter, we prove Theorems 3.2, 3.4, and 3.9 respectively.

\subsection{Theorem 3.2}
{\bf Proof of Theorem 3.2} Let the set of P-positions given in the theorem be \\
\quad $P_0=\{()\}, \quad P_1=\{(m,m,m) \ | \ m\geq 1\}$,\\
\quad $P_2=\{(m,M,m,M) \ | \ m<M\}$,\\
\quad $P_3=\{(M,m,a,b,m) \ | \ m<a\leq b<M\}$,\\
\quad $P_4=\{(m+1,M,m,m,M) \ |\ m\mbox{ is even}, m+2\leq M\}$,\\
\quad $P_5=\{(m+1,M,m,m,M) \ | \ m\mbox{ is odd}, m+3\leq M\}$,\\
\quad $P_6=\{(m+2,m+1,m,m,m+1) \ | \ m\mbox{ is odd}\}$,\\
\quad $P_7=\{(m+1,m+1,m,m,m+2) \ | \ m\mbox{ is odd}\}$,\\
\quad $P=P_0\cup P_1 \cup P_2\cup P_3\cup P_4\cup P_5\cup P_6\cup P_7$\\
 (adding the condition $a\leq b$ to $P_3$ due to symmetry). 
The conditions $m\neq a$ and $b\neq M$ imposed in $P_3$ ensure that a move cannot lead to a  configuration in $P_2$. 
 It can be proven that no move from a  configuration in $P_i$ can lead to a  configuration in 
$P_j$ for $0\leq i,j\leq 7$ individually. Since none of them are difficult and detailing them would be redundant, this paper omits their description. 
Below, we prove that if $(m,a,b,c,d)\not\in P$ and $m=\min\{m,a,b,c,d\}$, then there exists a move from $(m,a,b,c,d)$ to a P-position.

First, assume $m=a$. If $m=a=b$ or $m=a=d$, then there is a move to $(m,m,m)\in P_1$.
In other cases, we can assume without loss of generality that $m=a<b\leq d$. If $c=m$,  there is a move to $(m,b,m,b)\in P_2$.
If  $c>m$,
\begin{itemize}
\item If $b\geq m+2$ and $m$ is even, there is a move to $(m,m,b,m+1,b)\in P_4$.
\item If $b=m+1$ and $m$ is even, there is a move to $(m,m,m+1,m-1,m-1)\in P_7$.
\item If $b\geq m+3$ and $m$ is odd, there is a move to $(m,m,b,m+1,b)\in P_5$.
\item If $b=m+2$ and $m$ is odd, there is a move to $(m,m,m+2,m+1,m+1)\in P_7$.
\item If $b=m+1$ and $m$ is odd, there is a move to  $(m+1,m+1,m+1)\in P_1$ when $c=d=m+1$;  in other cases, there is a move to $(m,m,m+1,m+1,m+2)\in P_7$ or $(m,m,m+1,m+2,m+1)\in P_6$.
\end{itemize}

Next, consider the case where $m<a$. Without loss of generality, assume $a\leq d$. If $a=d$, then there is a move to $
(m,a,m,a)\in P_2$. If $a<d$, 
\begin{itemize}
\item If $m+d\leq a+b$, then since $m<m+(d-a)=d-(a-m)<d$, there is a move to $(m,a,m+d-a,m,d)\in P_3$.
\item If $m+d>a+b$ and $m<b$, then since $a<a+(b-m)$ and $b<b+(a-m)$, there is a move to $(m,a,b,m,a+b-m)\in P_3$.
\item If $m=b$, there is a move to $(m,a,m,a)\in P_2$. \quad$\Box$
\end{itemize}

\subsection{Theorem 3.4}
{\bf Proof of Theorem 3.4} Let the set of P-positions given in the theorem be \\
\quad $P_0=\{()\}, \quad P_1=\{(m,m,m,m) \ | \ m\geq 1\}$, \\
\quad $P_2=\{(1,M,a,b,M) \ | \  1\leq a<b\leq M, 1+M=a+b\}$,\\
\quad $P_3=\{(2,2p,p+1,p,2p-1) \ | \ p\geq 2\}$,\\
\quad $P=P_0\cup P_1 \cup P_2\cup P_3$.

The proof that no move from $(a,b,c,d,e)\in P$ can lead to any P-position is ommitted in this paper, as it is lengthy but not difficult. 

Below, we prove that if  $(a,b,c,d,e)\not\in P$, then there exists a move from $(a,b,c,d,e)$ to a P-position. Assume that $b=\max\{a,b,c,d,e\}$ and, without loss of generality, that $d\leq e$. If $e=1$, then $d=1$, so there is a move to $(1,1,1,1)\in P_1$. If $c=d$, there is a move to $(d,d,d,d)\in P_1$. Thus, we assume $e\geq 2$ and $c\neq d$ below. We proceed by case analysis based on the comparison of $1+e$ and $c+d$.
\begin{itemize}
\item When $1+e\geq c+d$, there is a move to $(1,c+d-1,c,d,c+d-1)\in P_2$.
\item When $1+e<c+d$ and $e\neq 2d-1$, there is a move to $(1,e,e-d+1,d,e)\in P_2$. Here, $e\neq 2d-1$ ensures that $e-d+1\neq d$.
\item When $1+e<c+d$ and $e=2d-1$, these conditions imply $c>d$. If $a=1$, then $(1,b,c,d,2d-1)$ can be moved to $(1,2d-1,d+1,d-1,2d-1)\in P_2$. If $a\geq 2$ and $b>e=2d-1$, then $(a,b,c,d,2d-1)$ can be moved to $(2,2d,d+1,d,2d-1)\in P_3$. If $a\geq 2$ and $b=e=2d-1$, then then $(a,2d-1,c,d,2d-1)$ can be moved to $(1,2d-1,c,2d-c,2d-1)\in P_2$. Here, we use the maximality of $b$ to ensure $c\leq b=2d-1$. \quad $\Box$
\end{itemize}

\subsection{Theorem 3.9}
{\bf Proof of Theorem 3.9} Let the set of P-positions given in the theorem be $P$, that is:\\
\quad $P_0=\{()\}, \quad P_1=\{(m,m,m,m,m,m,m)\ | \ m\geq 1\}$, \\
\quad $\tilde{P}_2=\{(1,M,a,M-a+1,\alpha,M-b+1,b,M) \ | \ (\ast)\}$,\\
\quad $N_2=\{(1,2p-1,p,p,2p-1,p,p,2p-1) \ | \ p\geq 1\}$,\\
\quad $P_2=\tilde{P}_2\setminus N_2$, \quad $P=P_0\cup P_1\cup P_2$, \\
where the condition $(\ast)$ be defined as $a,b\leq M$ and $\alpha=\min\{M,a+b-1\}$.

Let $\Omega=\{(a,b,c,d,e,f,g,h) \ | \ a,b,c,d,e,f,g,h\geq 1\}$, which represents the set of configurations for ${\rm SCN}(8,6)$ where all eight piles are remaining. We give a proof of this theorem in four steps:
\begin{enumerate}
\renewcommand{\theenumi}{(\roman{enumi})}
\item $\omega\in\tilde{P}_2$ can be moved to some configuration in $P_1$ if and only if $\omega\in N_2$.
\item There is no move from one configuration in $\tilde{P}_2$ to another configuration in $\tilde{P}_2$.
\item $\omega\in\Omega\setminus\tilde{P}_2$ can be moved to some configuration in $\tilde{P}_2$.
\item If  $\omega\in\Omega\setminus\tilde{P}_2$ can be moved to some configuration in $N_2$, then $\omega$ can also be moved to some configuration in $P_1\cup P_2$.
\end{enumerate}
(i)  $(1,2p-1,p,p,2p-1,p,p,2p-1)\in N_2$ can be moved to $(p,p,p,p,p,p,p)\in P_1$. Hereafter, we assume $\omega=(1,M,a,M-a+1,\alpha,M-b+1,b,M)\in \tilde{P}_2$ can be moved to $(m,m,m,m,m,m,m)\in P_1$ and prove that $\omega\in N_2$.  We note that all but one of the piles have at least $m$ stones in $\omega$. Thus, $m\leq a,b, M-a+1$ and $M-b+1$.  
\begin{itemize}
\item If $m=1=M$, then $a=b=1$ and $\omega=(1,1,1,1,1,1,1,1)\in N_2$.
\item If $m=M=a\leq M-a+1$, then $a=1$ that implies $\omega=(1,1,1,1,1,1,1,$\\$1)\in N_2$. 
\item If $m=a=M-a+1$, then $a\leq b$ and $M-a+1\leq M-b+1$ that implies $a=b$, $M=2a-1$ and $\omega=(1,2a-1,a,a,2a-1,a,a,2a-1)\in N_2$.
\item Suppose $m=M-a+1=\min\{M,a+b-1\}$. If $M-a+1=M\leq M-b+1$ then $a=b=1$, that implies  $\omega=(1,1,1,1,1,1,1,1)\in N_2$.  If $M-a+1=a+b-1\leq M$, then $M-a+1\leq M-b+1$ and $a+b-1\leq b$ implies $a=b=M=1$ and $\omega=(1,1,1,1,1,1,1,1)\in N_2$.
\end{itemize} 
(ii) Let $\Omega'=\{(a,b,c,d,e,f,g,h) \ | \ a,b,c,d,e,f,g,h\geq 0\}$, which represents the set of configurations of ${\rm CN}(8,6)$.
We define a bijective mapping $\tau$ from $\Omega$ to $\Omega'$, that maps $(a,b,c,d,e,f,g,h)\in\Omega$ to $(a-1,b-1,c-1,d-1,e-1,f-1,g-1,h-1)\in\Omega'$.
 It is important to observe that $\tau(\tilde{P}_2)$ coincides precisely with the set of P-positions of  ${\rm CN}(8,6)$.  
If $\omega\in\tilde{P}_2$ can be moved to another configuration in $\tilde{P}_2$, 
then $\tau(\omega)$ in $\tau(\tilde{P}_2)$ can be moved to another P-position of ${\rm CN}(8,6)$, that is a contradiction.\\
(iii) For any $\omega\in\Omega\setminus\tilde{P}_2$, $\tau(\omega)\in\Omega'$ is a N-position of ${\rm CN}(8,6)$, so $\tau(\omega)$ can be moved to some P-position $\omega'$ of  ${\rm CN}(8,6)$. Then $\omega$ can be moved to $\tau^{-1}(\omega')\in \tilde{P}_2$. \\
(iv) We will prove that if $\omega\in \Omega\setminus \tilde{P}_2$ can be moved to $(1,2p-1,p,p,2p-1,p,p,2p-1)\in N_2$ (where $p\geq 1$) then $\omega$ can also be moved to another  configuration in $P$. Below, we classify the original configurations in $\Omega\setminus \tilde{P}_2$ that can be moved to $(1,2p-1,p,p,2p-1,p,p,2p-1)$ into five cases.
 
\begin{itemize}
\item When the resulting  configuration corresponds to $p=1$, that is $(1,1,1,1,$\\$1,1,1,1)$, the original  configuration must contain two consecutive piles with exactly one stone each. Therefore, it is also possible to move from the original  configuration to $(1,1,1,1,1,1,1)\in P_1$.
\end{itemize}

Based on the above, in the following four cases, we also assume $p\geq 2$ (which implies $p<2p-1$).
\begin{itemize}
\item Suppose $(1,2p-1,C,D,E,F,G,H)$ can be moved to $(1,2p-1,p,p,2p-1,p,p,2p-1)$. In this case 
\begin{itemize}
\item If $C=D=p$, the original  configuration can also be moved to $(p,p,p,p,p,p,p)\in P_1$.
\item If $C>p$, it can be moved to $(1,2p-1,p+1,p-1,2p-1,p,p,2p-1)\in P_2$.
\item If $D>p$, it can be moved to $(1,2p-1,p-1,p+1,2p-2,p,p,2p-1)\in P_2$.
\end{itemize}
\item Now, suppose $(A,2p-1,p,D,E,F,G,H)$ can be moved to $(1,2p-1,p,p,2p-1,p,p,2p-1)$. In this case 
\begin{itemize}
\item If $F=G=p$, the original  configuration can also be moved to $(p,p,p,p,p,p,p)\in P_1$.
\item If $F>p$, it can be moved to $(1,2p-1,p,p,2p-2,p+1,p-1,2p-1)\in P_2$.
\item If $G>p$, it can be moved to $(1,2p-1,p,p,2p-1,p-1,p+1,2p-1)\in P_2$. 
\end{itemize}
\item Next, suppose $(A,B,p,p,E,F,G,H)$ can be moved to $(1,2p-1,p,p,2p$\\$-1,p,p,2p-1)$. In this case, the original  configuration can also be moved to $(p,p,p,p,p,p,p)\in P_1$. 
\item Lastly, suppose $(A,B,C,p,2p-1,F,G,H)$ can be moved to $(1,2p-1,p,p,2p-1,p,p,2p-1)$. 
\begin{itemize}
\item If $C=p$, it can also be moved to $(p,p,p,p,p,p,p)\in P_1$. 
\item If $G>p$, it can be moved to $(1,2p-1,p,p,2p-1,p-1,p+1,2p-1)\in P_2$.
\item If $F=G=p$, it can be moved to $(p,p,p,p,p,p,p)\in P_1$.
\end{itemize}

\ So we assume below that $C>p$, $G=p$ and $F>p$. Under these assumptions, if $H=2p-1$, it can be moved to $(1,2p-1,p+1,p-1,2p-1,p,p,2p-1)\in P_2$. So we also assume $H>2p-1$.

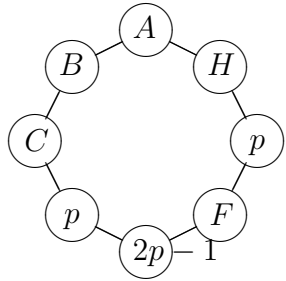
\begin{figure}[ht]
\unitlength.1pt
\begin{picture}(1100,1200)(-1400,-150)
\put(500,80){\circle{200}}
\put(220,220){\circle{200}}
\put(220,780){\circle{200}}
\put(80,500){\circle{200}}
\put(500,920){\circle{200}}
\put(780,220){\circle{200}}
\put(780,780){\circle{200}}
\put(920,500){\circle{200}}
\put(309,175){\line(2,-1){102}}
\put(691,175){\line(-2,-1){102}}
\put(309,825){\line(2,1){102}}
\put(691,825){\line(-2,1){102}}
\put(175,309){\line(-1,2){54}}
\put(825,309){\line(1,2){54}}
\put(175,691){\line(-1,-2){54}}
\put(825,691){\line(1,-2){54}}
\put(450,890){\large{$A$}}
\put(450,50){\large{$2p-1$}}
\put(190,190){\large{$p$}}
\put(740,180){\large{$F$}}
\put(40,460){\large{$C$}}
\put(890,470){\large{$p$}}
\put(170,750){\large{$B$}}
\put(730,750){\large{$H$}}
\put(-550,-150){$A\geq 1,B\geq 2p-1,C>p,F>p,H>2p-1$}
\end{picture}
\caption{ The intermediate stage of the proof  ($\ast\ast$)}
\end{figure}

\ We now list any remaining cases ($\ast\ast$) that have not yet been proven as Figure 12. Here, if we assume $F\leq 2p-1$:
\begin{itemize}
\item If $F\leq C$, then the original  configuration can also be moved to $(1,2p-1,F,2p-F,2p-1,F,2p-F,2p-1)\in P_2$.
\item If $F>C$, when $B=2p-1$ it can be moved to  $(1,2p-1,C,2p-C,2p-1,p,p,2p-1)\in P_2$ and when $B>2p-1$ it can be moved to $(1,2p,p+1,p,2p-1,p+2,p-1,2p)\in P_2$.
\end{itemize}

\ Thus, it is sufficient to consider the case where $F>2p-1$. If we further assume $B>2p-1$, then the original  configuration can also be moved to $(1,2p,p+1,p,2p-1,p+2,p-1,2p)\in P_2$. If we assume $B=2p-1$ and $C\leq 2p-1$, then it can be moved to $(1,2p-1,C,2p-C,2p-1,p,p,2p-1)\in P_2$. Therefore, for the remaining proof, we can assume the situation  ($\ast\!\ast\!\ast$) illustrated in Figure 13.
\begin{figure}[ht]
\unitlength.1pt
\begin{picture}(1100,1200)(-1400,-150)
\put(500,80){\circle{200}}
\put(220,220){\circle{200}}
\put(220,780){\circle{200}}
\put(80,500){\circle{200}}
\put(500,920){\circle{200}}
\put(780,220){\circle{200}}
\put(780,780){\circle{200}}
\put(920,500){\circle{200}}
\put(309,175){\line(2,-1){102}}
\put(691,175){\line(-2,-1){102}}
\put(309,825){\line(2,1){102}}
\put(691,825){\line(-2,1){102}}
\put(175,309){\line(-1,2){54}}
\put(825,309){\line(1,2){54}}
\put(175,691){\line(-1,-2){54}}
\put(825,691){\line(1,-2){54}}
\put(450,890){\large{$A$}}
\put(450,50){\large{$2p-1$}}
\put(190,190){\large{$p$}}
\put(730,180){\large{$F$}}
\put(40,460){\large{$C$}}
\put(890,470){\large{$p$}}
\put(160,750){\large{$2p-1$}}
\put(730,750){\large{$H$}}
\put(-550,-150){$A\geq 1,C>2p-1,F>2p-1,H>2p-1$}
\end{picture}
\caption{ The intermediate stage of the proof  ($\ast\!\ast\!\ast$) }
\end{figure}
\end{itemize}

Under these assumptions, we prove that the original  configuration can be moved to a  configuration in $P_2$ by considering four cases.  Note that in each resulting  configuration shown below, the rightmost pile $G$ in the diagram will have one stone, instead of the pile $A$ as previously discussed. 
\begin{enumerate}
\item When $C,F,H\geq\min\{3p-2,2p+A-2\}$:
\begin{itemize}
\item If $3p-2\leq 2p+A-2$, it can be moved to $(p,2p-1,3p-2,p,2p-1,3p-2,1,3p-2)\in P_2$.
\item If $3p-2>2p+A-2$, it can be moved to $(A,2p-1,2p+A-2,A,2p-1,2p+A-2,1,2p+A-2)\in P_2$. Here $A\neq 2p-1$ is confirmed by $A<p<2p-1$.
\end{itemize}
\item When $C\leq F,H$ and $C<3p-2$ and $C<2p+A-2$, it can be moved to $(C-2p+2,2p-1,C,C-2p+2,2p-1,C,1,C)\in P_2$. Here $2p-1\neq C-2p+2$ is confirmed by $C<3p-2<4p-3$. 
\item When $F\leq C,H$ and $F<3p-2$ and $F<2p+A-2$, it can be moved to $(F-2p+2,2p-1,F,F-2p+2,2p-1,F,1,F)\in P_2$.
Here, $2p-1\neq F-2p+2$ is confirmed by $F<3p-2<4p-3$.
\item When $H\leq C,F$ and $H<3p-2$ and $H<2p+A-2$, we consider two subcases based on the on the comparison of $p$ and $A$:
\begin{itemize}
\item When $p\geq A$, it can be moved to $(A,H-A+1,H,A,H-A+1,H,1,H)\in P_2$. Here $H-A+1\neq A$ is confirmed by $H>2p-1\geq 2A-1$.
\item When $p<A$, it can be moved to $(A,H-A+1,H,p,H-p+1,H,1,H)\in P_2$. Here $p\neq H-p+1$ is confirmed by $H>2p-1$, and $\min\{H,A+(H-p+1)-1\}=H$ is confirmed by $A+(H-p+1)-1=H+(A-p)>H$. \quad$\Box$
\end{itemize}
\end{enumerate}

\section{Concluding Remarks}

\subsection{Difference between Circular Nim and Shrinking Circular Nim}

We demonstrated the winning conditions for several cases of Shrinking Circular Nim. In each case, while the winning conditions for Circular Nim ${\rm CN}(n, k)$ are given by a simple formula, it was found that for the corresponding Shrinking Circular Nim ${\rm SCN}(n, k)$ with the same $(n, k)$, additional conditions are attached to the formula, and exceptional P-positions arise. This difference stems from the fact that in Shrinking Circular Nim, the circle may contract, which results in the winning conditions of ${\rm SCN}(n', k)$ ($n'<n$) being layered upon one another. Even if the winning conditions for Circular Nim are resolved for other values of $(n, k)$ in the future, it is anticipated that the resolution will become more difficult in Shrinking Circular Nim as the number of piles $n$ increases.

On the other hand, since the rule for removing stones in Shrinking Circular Nim is the same as in Circular Nim, it is also speculated that the essential aspect of the winning conditions will not change, even with an increase in exceptional configurations. One such example was the resolution of ${\rm SCN}(8, 6)$. In fact, many of the P-positions of Shrinking Circular Nim obtained in this paper resemble the P-positions of Circular Nim by subtracting one stone from each pile, leading to the optimistic observation that as the solution for ${\rm CN}(n, k)$ progresses, ${\rm SCN}(n, k)$ may also be resolved.

Neither Circular Nim nor Shrinking Circular Nim is inherently superior as a game, but it should be noted that the increased complexity of the winning conditions does not hinder the actual playing of the game. In fact, the complexity may make the game more enjoyable. Furthermore, when using this game for verification in machine learning or similar fields, the existence of exceptional  configurations is considered valuable.

\subsection{Conclusion}

This paper introduces a game called Shrinking Circular Nim, which is a variation of the Circular Nim game where some of the rules have been modified. It examines the winning conditions when the number of piles is six or fewer, as well as in special cases where the winning conditions for Circular Nim have already been resolved. As a result, we were able to fully obtain the set of P-positions for cases with five or fewer piles. This set of P-positions contains exceptional cases compared to those of Circular Nim, and while it may lack mathematical simplicity, this complexity contributes to the strategic depth and enjoyment of the game. Future challenges include obtaining the winning conditions for a more general number of piles and deepening the understanding of the general relationship between Circular Nim and Shrinking Circular Nim.

\section*{Acknowledgements}
This work was supported by JSPS KAKENHI Grant Number JP21K12191.

%
%
%
 \bibliographystyle{splncs04}
 \bibliography{SCN2024forarXiv}
\end{document}